%% file: preprint.tex
\documentclass[runningheads]{llncs}
\usepackage[breaklinks, colorlinks, linkcolor=blue, citecolor=blue, urlcolor=black]{hyperref}
\usepackage{color}

\usepackage{graphicx}
\usepackage{tabularx}
\usepackage{tabularcalc}
\usepackage{booktabs}
\usepackage{multirow}
\usepackage{algorithm}
\usepackage{algpseudocode}
\usepackage{tikz}
\usepackage{forest}
\usepackage{pgfplots}
\usepackage{xltabular}
\usetikzlibrary{arrows.meta}
\usepackage{amsmath}
\usepackage{amssymb}
\usepackage{xspace}
\usepackage{xcolor}
\usepackage{varioref}
\usepackage{cleveref}
\usepackage{aligned-overset}
\usepackage[shortlabels,inline]{enumitem}
\usepackage{boxedminipage2e}
\usepackage{makecell}

\definecolor{pbblue}{RGB}{68,119,170}
\definecolor{pbcyan}{RGB}{102,204,238}
\definecolor{pbgreen}{RGB}{34,136,51}
\definecolor{pbred}{RGB}{238,102,119}

\input{macros-general.tex}

\input{macros-proofcplx.tex}


\newcommand{\RNum}[1]{\uppercase\expandafter{\romannumeral #1\relax}}
\newcommand{\solver}[1]{\textsc{#1}\xspace}

\newcommand{\scipversion}{9.3}
\newcommand{\soplexversion}{8.0}

\newcommand{\sciphash}{(hash: \solver{b7906251d5})}
\newcommand{\soplexhash}{(hash: \solver{7b9bd461})}

\newcommand{\exactsciphash}{(hash: \solver{c7bd4be7b9})}

\newcommand{\scip}{\solver{SCIP}}
\newcommand{\soplex}{\solver{SoPlex}}

\newcommand{\scipv}{\scip~\scipversion\xspace}
\newcommand{\soplexv}{\soplex~\soplexversion\xspace}

\newcommand{\exactscip}{\solver{Exact-SCIP}}
\newcommand{\exactsoplex}{\solver{Exact-SoPlex}}

\newcommand{\fpeasy}{\solver{FPEasy}}
\newcommand{\numdiff}{\solver{NumDiff}}

\newcommand{\bestsol}{\solver{BESTSOLUTION}}
\newcommand{\nodedel}{\solver{NODEDELETE}}
\newcommand{\nodefeas}{\solver{NODEFEASIBLE}}
\newcommand{\nodeinfeas}{\solver{NODEINFEASIBLE}}
\newcommand{\bandb}{branch-and-bound\xspace}

\newcommand{\fpsolver}{floating-point solver\xspace}

\newcommand{\bernd}{\solver{bernd2}}
\newcommand{\dfnload}{\solver{dfn6\_load}}

\newcommand{\myorcidlink}[1]{\,\href{https://orcid.org/#1}{\raisebox{-0.45ex}{\includegraphics[width=1.8ex]{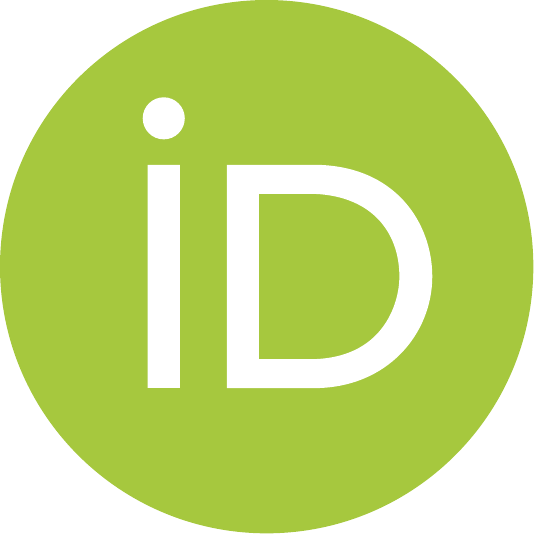}}}}

\definecolor{pbblue}{RGB}{68,119,170}
\definecolor{pbcyan}{RGB}{102,204,238}
\definecolor{pbgreen}{RGB}{34,136,51}
\definecolor{pbred}{RGB}{238,102,119}

\usepackage{todonotes,xspace}

\makeatletter
\def\orcidID#1{\href{http://orcid.org/#1}{\protect\raisebox{-1.25pt}{\protect\includegraphics{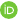}}}}
\makeatother

\begin{document}

\title{Analyzing the numerical correctness of \bandb decisions for mixed-integer programming}
\author{
    Alexander Hoen\inst{1}\orcidID{0000-0003-1065-1651} \and
    Ambros Gleixner\inst{1,2}\orcidID{0000-0003-0391-5903} 
    \institute{%
      Zuse Institute Berlin, Takustr.~7, 14195 Berlin, Germany\\
      \email{hoen@zib.de}
      \and
      HTW Berlin, 10313 Berlin, Germany\\
      \email{gleixner@htw-berlin.de} \\
    }
}

\authorrunning{A. Hoen, A. Gleixner}
\titlerunning{Analyzing the numerical correctness of \bandb decisions for MIP} 

\maketitle

\begin{abstract}
  \input{abstract}

\end{abstract}

\input{sections/introduction}

\input{sections/basics}
\input{sections/experiments}

\input{sections/conclusion}

\vspace{0.5cm}
\noindent{\small
\input{sections/acknowledgements}

}

\bibliographystyle{splncs04}
\bibliography{bibliography}

\end{document}

%% file: macros-general.tex
\usepackage{ifthen}
\input{doctypedetection.tex}

\usepackage{xspace}
\ifthenelse
{\boolean{detectedElsevier} \or \boolean{detectedSIAM} 
  \or \boolean{detectedLIPIcs}}
{}
{\usepackage{varioref}}

\makeatletter
\AtBeginDocument{%
  \def\doi#1{\url{https://doi.org/#1}}}
\makeatother

\DeclareMathAlphabet{\mathsfsl}{OT1}{cmss}{m}{sl}

\DeclareRobustCommand{\BibTeX}{%
  {\normalfont B\kern-.05em{\scshape i\kern-.025em b}\kern-.08em \TeX}%
}

\ifthenelse{\boolean{detectedToC}}{}
{

}

\newcommand{\MAXOFEXPR}[2][]{\max_{#1} \left\{ #2 \right\}}
\newcommand{\MINOFEXPR}[2][]{\min_{#1} \left\{ #2 \right\}}
\newcommand{\Maxofexpr}[2][]{\max_{#1} \bigl\{ #2 \bigr\}}
\newcommand{\Minofexpr}[2][]{\min_{#1} \bigl\{ #2 \bigr\}}

\newcommand{\MAXOFSET}[3][:]%
     {\ifthenelse{\equal{#1}{;}}%
     {\MAXOFEXPR{ #2 \,;\, #3 }}
     {\ifthenelse{\equal{#1}{:}}%
     {\MAXOFEXPR{ #2 \,:\, #3 }}
     {\max \twincommandJN{\left\{}{#2}{\left#1}{\right}{\,#3}{\right\}}}}}
\newcommand{\MINOFSET}[3][:]%
     {\ifthenelse{\equal{#1}{;}}%
     {\MINOFEXPR{ #2 \,;\, #3 }}
     {\ifthenelse{\equal{#1}{:}}%
     {\MINOFEXPR{ #2 \,:\, #3 }}
     {\min \twincommandJN{\left\{}{#2}{\left#1}{\right}{\,#3}{\right\}}}}}

\newcommand{\Maxofset}[3][:]%
     {\ifthenelse{\equal{#1}{;}}%
     {\Maxofexpr{ #2 \,;\, #3 }}
     {\ifthenelse{\equal{#1}{:}}%
     {\Maxofexpr{ #2 \,:\, #3 }}
     {\max \twincommandJN{\bigl\{}{#2}{\bigl#1}{\bigr}{\,#3}{\bigr\}}}}}
\newcommand{\Minofset}[3][:]%
     {\ifthenelse{\equal{#1}{;}}%
     {\Minofexpr{ #2 \,;\, #3 }}
     {\ifthenelse{\equal{#1}{:}}%
     {\Minofexpr{ #2 \,:\, #3 }}
     {\min \twincommandJN{\bigl\{}{#2}{\bigl#1}{\bigr}{\,#3}{\bigr\}}}}}

\DeclareMathOperator{\Expop}{E}

\ifthenelse{\boolean{detectedLMCS}}
{}
{}

\newcommand{\twincommandJN}[6]%
    {#1#2#3\vphantom{#2#5}\mspace{-2.05mu}#4.#5#6}

\newcommand{\CondExp}[2]%
    {\Expop\twincommandJN{\bigl[}{#1}{\bigl|}{\bigr}{\,#2}{\bigr]}}
\newcommand{\CONDEXP}[2]%
     {\Expop\twincommandJN{\left[}{#1}{\left|}{\right}{\,#2}{\right]}}

\newcommand{\Condprob}[3][]%
    {\Pr_{#1}\twincommandJN{\bigl[}{#2}{\bigl|}{\bigr}{\,#3}{\bigr]}}
\newcommand{\CONDPROB}[3][]%
    {\Pr_{#1}\twincommandJN{\left[}{#2}{\left|}{\right}{\,#3}{\right]}}

\newcommand{\Set}[1]{\bigl\{ #1 \bigr\}}

\newcommand{\Setdescr}[3][|]%
     {\ifthenelse{\equal{#1}{;}}%
     {\Set{ #2 \,;\, #3 }}
     {\ifthenelse{\equal{#1}{:}}%
     {\Set{ #2 \,:\, #3 }}
     {\twincommandJN{\bigl\{}{#2\,}{\bigl#1}{\bigr}{\,#3}{\bigr\}}}}}
\newcommand{\SETDESCR}[3][|]%
     {\twincommandJN{\left\{}{#2\,}{\left#1}{\right}{\,#3}{\right\}}}

\newcommand{\Setdescrbrackets}[3][|]%
     {\twincommandJN{\bigl[}{#2}{\bigl#1}{\bigr}{\,#3}{\bigr]}}
\newcommand{\SETDESCRBRACKETS}[3][|]%
     {\twincommandJN{\left[}{#2}{\left#1}{\right}{\,#3}{\right]}}

\newcommand{\nvar}{n}

\newcommand{\nclause}{m}

\newcommand{\clwidth}{k}

\newcommand{\randkcnfnclwrepl}[3][\clwidth]%
        {\ensuremath{\mathcal{F}^{#2, #3}_{#1}}}
\newcommand{\randkcnfnclwreplstd}%
        {\randkcnfnclwrepl{\clwidth}{\nvar}{\nclause}}

\newcommand{\complclassformat}[1]%
        {\textrm{\upshape{\textsf{#1}}}\xspace}
\newcommand{\cocomplclass}[1]%
        {\textrm{\upshape{\textsf{co#1}}}\xspace}

\newcommand{\DTIMEadviceclass}[2]%
    {\ensuremath{\complclassformat{DTIME}\bigl(#1\bigr)/{#2}}}

\newcommand{\PCPalph}[5]%
    {\ensuremath{\complclassformat{PCP}_{{#1},{#2}}[{#3}, {#4}, {#5}]}}
\newcommand{\PCP}[4]%
    {\ensuremath{\complclassformat{PCP}_{{#1},{#2}}[{#3}, {#4}]}}

\ifthenelse{\boolean{detectedToC}}{}
  { 
    }
\ifthenelse{\boolean{detectedToC}}{}{

}

\ifthenelse{\boolean{detectedLIPIcs} \or \boolean{detectedIJCAI}}
{\renewcommand{\st}{\errmessage{Please do not use st}}}
{\ifthenelse{\isundefined{\st}}
  {\newcommand{\st}{such that\xspace}}
  {}}      

\ifthenelse{\boolean{detectedIEEE}}{}{}

\ifthenelse
{\isundefined{\refeq}}
{\newcommand{\refeq}[1]{\eqref{#1}}}
{\renewcommand{\refeq}[1]{\eqref{#1}}}

%% file: doctypedetection.tex
\usepackage{ifthen}

\provideboolean{detectedSTOC}
\provideboolean{detectedFOCS}
\provideboolean{detectedElsevier}
\provideboolean{detectedNOW}
\provideboolean{detectedLMCS}
\provideboolean{detectedIEEE}
\provideboolean{detectedPoster}
\provideboolean{detectedSIAM}
\provideboolean{detectedLNCS}
\provideboolean{detectedACM}
\provideboolean{detectedACMconf}
\provideboolean{detectedSigplanconf}
\provideboolean{detectedToC}
\provideboolean{detectedLIPIcs}
\provideboolean{detectedAAAI}
\provideboolean{detectedIJCAI}
\provideboolean{detectedCompCplx}
\provideboolean{detectedEasyChair}
\provideboolean{detectedArticle}
\provideboolean{detectedReport}
\provideboolean{detectedThesis}

\makeatletter

\@ifclassloaded{sig-alternate}
{\setboolean{detectedSTOC}{true}}
{\setboolean{detectedSTOC}{false}}

\@ifclassloaded{elsarticle}
{\setboolean{detectedElsevier}{true}}
{\setboolean{detectedElsevier}{false}}

\@ifclassloaded{now}
{\setboolean{detectedNOW}{true}}
{\setboolean{detectedNOW}{false}}

\@ifclassloaded{lmcs}
{\setboolean{detectedLMCS}{true}}
{\setboolean{detectedLMCS}{false}}

\@ifclassloaded{IEEEtran} {
  \setboolean{detectedIEEE}{true}
  \ifCLASSOPTIONconference {
    \setboolean{detectedFOCS}{true}
  }
  \else {
    \setboolean{detectedFOCS}{false}
  }
  \fi
}
{
  \setboolean{detectedFOCS}{false}
  \setboolean{detectedIEEE}{false}
}

\@ifclassloaded{siamart171218}
{\setboolean{detectedSIAM}{true}}
{\setboolean{detectedSIAM}{false}}

\@ifclassloaded{llncs}
{\setboolean{detectedLNCS}{true}}
{\setboolean{detectedLNCS}{false}}

\@ifclassloaded{acmsmall}
{\setboolean{detectedACM}{true}}
{\setboolean{detectedACM}{false}}

\@ifclassloaded{acmart}
{\setboolean{detectedACMconf}{true}}
{\setboolean{detectedACMconf}{false}}

\@ifclassloaded{sigplanconf}
{\setboolean{detectedSigplanconf}{true}}
{\setboolean{detectedSigplanconf}{false}}

\@ifclassloaded{toc}
{\setboolean{detectedToC}{true}}
{\setboolean{detectedToC}{false}}

\@ifclassloaded{lipics}
{\setboolean{detectedLIPIcs}{true}}
{\@ifclassloaded{lipics-v2019}
  {\setboolean{detectedLIPIcs}{true}}
  {\@ifclassloaded{oasics-v2019}
    {\setboolean{detectedLIPIcs}{true}}
    {\setboolean{detectedLIPIcs}{false}}
  }
}

\@ifclassloaded{cc}
{\setboolean{detectedCompCplx}{true}}
{\setboolean{detectedCompCplx}{false}}

\@ifclassloaded{easychair}
{\setboolean{detectedEasyChair}{true}}
{\setboolean{detectedEasyChair}{false}}

\@ifpackageloaded{aaai}
{\setboolean{detectedAAAI}{true}}        
{\@ifpackageloaded{aaai18}
  {\setboolean{detectedAAAI}{true}}
  {\@ifpackageloaded{aaai20}       
    {\setboolean{detectedAAAI}{true}}
    {\setboolean{detectedAAAI}{false}}}}

\@ifpackageloaded{ijcai18}
{\setboolean{detectedIJCAI}{true}}        
{\@ifpackageloaded{ijcai19}
  {\setboolean{detectedIJCAI}{true}}        
  {\setboolean{detectedIJCAI}{false}}}

\@ifclassloaded{sciposter}
{\setboolean{detectedPoster}{true}}
{\setboolean{detectedPoster}{false}}

\@ifclassloaded{article}
{\setboolean{detectedArticle}{true}}
{\setboolean{detectedArticle}{false}}

\makeatother

\ifthenelse{\not \isundefined{\examen} 
  \and \not \isundefined{\disputationsdatum} 
  \and \not \isundefined{\disputationslokal}}   
  {\setboolean{detectedThesis}{true}}
  {\setboolean{detectedThesis}{false}}

\ifthenelse{\boolean{detectedArticle} \or \boolean{detectedThesis}
  \or \boolean{detectedSTOC} \or \boolean{detectedFOCS}
  \or \boolean{detectedSIAM} \or \boolean{detectedIEEE}
  \or \boolean{detectedPoster}}
{\setboolean{detectedReport}{false}}
{\setboolean{detectedReport}{true}}

%% file: macros-proofcplx.tex
\ifthenelse{\boolean{detectedACMconf}}
{}
{}

\ifthenelse{\boolean{detectedACMconf}}
{
 }
{
 }

\newcommand{\SETSOFVARSORLIT}[2]%
        {\mathit{#1}\left({#2}\right)}
\newcommand{\setsofvarsorlit}[2]%
        {\mathit{#1}({#2})}
\newcommand{\Setsofvarsorlit}[2]%
        {\mathit{#1}\bigl({#2}\bigr)}

\newcommand{\derivabbrev}[2]{\bigl( #1 \vdash #2 \bigr)}
\newcommand{\derivabbrevsmall}[2]{( #1 \vdash #2 )}
\newcommand{\derivabbrevcompact}[2]{\bigl( #1 \vdash #2 \bigr)}

\newcommand{\refutabbrevsmall}[1]{\derivabbrevsmall{#1}{\!\bot}}
\newcommand{\refutabbrevcompact}[1]{\derivabbrevcompact{#1}{\!\bot}}

\newcommand{\genericrefsmall}[3]%
    {{\mathit{#1}}_{#2}\refutabbrevsmall{#3}}
\newcommand{\genericrefcompact}[3]%
    {{\mathit{#1}}_{#2}\refutabbrevcompact{#3}}
\newcommand{\genericderiv}[4]%
    {{\mathit{#1}}_{#2}\derivabbrev{#3}{#4}}
\newcommand{\genericderivsmall}[4]%
    {{\mathit{#1}}_{#2}\derivabbrevsmall{#3}{#4}}
\newcommand{\genericderivcompact}[4]%
    {{\mathit{#1}}_{#2}\derivabbrevcompact{#3}{#4}}
\newcommand{\generictaut}[3]%
    {{\mathit{#1}}_{#2}\derivabbrev{}{#3}}
\newcommand{\generictautcompact}[3]%
    {{\mathit{#1}}_{#2}\derivabbrevcompact{}{#3}}
\newcommand{\generictautsmall}[3]%
    {{\mathit{#1}}_{#2}\derivabbrevsmall{}{#3}}

\newcommand{\formulaformat}[1]{\mathit{#1}}

\newcommand{\extendedversion}[1]{\widetilde{#1}}

\newcommand{\epopnot}[1]%
    {\extendedversion{\formulaformat{POP}}_{#1}}

\newcommand{\elopnot}[1]%
    {\extendedversion{\formulaformat{LOP}}_{#1}}

\newcommand{\ephpnot}[2]%
    {\vphantom{\extendedversion{\formulaformat{PHP}}}
      {\smash{\extendedversion{\formulaformat{PHP}}}
        \vphantom{\formulaformat{PHP}}}^{#1}_{#2}}

\newcommand{\efphpnot}[2]%
    {\vphantom{\extendedversion{\formulaformat{FPHP}}}
      {\smash{\extendedversion{\formulaformat{FPHP}}}
        \vphantom{\formulaformat{FPHP}}}^{#1}_{#2}}

\newcommand{\ontophpnot}[2]%
    {\formulaformat{Onto}\text{-}\formulaformat{PHP}^{#1}_{#2}}
\newcommand{\ontofphpnot}[2]%
    {\formulaformat{Onto}\text{-}\formulaformat{FPHP}^{#1}_{#2}}

\newcommand{\graphontophpnot}[1][G]%
    {\text{$\formulaformat{Onto}$-$\formulaformat{PHP}$}({#1})}

\newcommand{\perfectmatchingnot}[1][G]%
    {\formulaformat{PM}({#1})}

%% file: abstract.tex
Most state-of-the-art branch-and-bound solvers for mixed-integer linear programming rely on limited-precision floating-point arithmetic and use numerical tolerances when reasoning about feasibility and optimality during their search.
While the practical success of floating-point MIP solvers bears witness to their overall numerical robustness, it is well-known that numerically challenging input can lead them to produce incorrect results.
Even when their final answer is correct, one critical question remains: 
Were the individual decisions taken during branch-and-bound justified, i.e., can they be verified in exact arithmetic?
In this paper, we attempt a first such \emph{a posteriori} analysis of a pure LP-based branch-and-bound solver by checking all intermediate decisions critical to the correctness of the result: accepting solutions as integer feasible, declaring the LP relaxation infeasible, and pruning subtrees as suboptimal.
Our computational study in the academic MIP solver \scip confirms the expectation that in the overwhelming majority of cases, all decisions are correct.
When errors do occur on numerically challenging instances, they typically affect only a small, typically single-digit, amount of leaf nodes that would require further processing.

%% file: sections/introduction.tex
\section{Introduction}
\label{sec::introdcution}

Mixed Integer Programming (MIP) solvers have evolved into highly sophisticated software capable of solving a growing number of large-scale and complex problems efficiently.
They typically implement a branch-and-bound scheme based on a linear programming (LP) relaxation.  
For performance reasons, nearly all solvers utilize floating-point arithmetic in order to speed up computations and control memory requirements for representing values of rational numbers.
These computational advantages of floating-point arithmetic, however, also imply limitations in precision: not all rational numbers can be represented exactly, leading to tiny round-off errors that can accumulate in magnitude over the course of the algorithm.
Therefore, solvers introduce tolerances for feasibility and a zero tolerance for deciding equality between two numbers. 
Additionally, by default, commercial solvers operate with relative gaps between $0.01\%$ and $0.0001\%$.
%
%

It is well documented that sometimes accumulated rounding errors can cause solvers to incorrectly declare infeasibility, feasibility, or optimality for suboptimal solutions~\cite{Cooketal11ExactRationalMIPSolver,Klotz2014,PaxianBiere-POS23,Steffy2011}.
Further, the MIPLIB~2017 instance collection labels 193 instances with problematic numerics with a tag ``numerics'', partly because different solvers reported inconsistent results during tests for the compilation process~\cite{GleixnerHendelGamrathetal2021MIPLIB2017}.
%
%
As a result, for applications that require absolute certainty, such as chip design verification \cite{Achterberg07Thesis}, combinatorial auctions \cite{VriesVohra03CombinatorailAuctions}, or computational proofs in experimental mathematics~\cite{EiflerGleixnerPulaj2022}, the level of trust provided by floating-point solvers is often not sufficient.

Nevertheless, it is widely accepted that the way that mature floating-point MIP solvers handle numerics is overall very robust.
This judgment is not only due to the fact that for the vast majority of industrial MIP applications, it may be acceptable to work with solutions that are slightly suboptimal or slightly infeasible.
The wide-spread use of floating-point MIP solvers and the fact that the general methodology of handling round-off errors in MIP solvers has not changed substantially over many years of development rather suggests that in the vast majority of cases, their answers may indeed be correct in an exact sense. 

However, to the best of our knowledge, the literature holds no record of a systematic computational investigation of this fundamental methodological question: To what extent do approximate floating-point solvers provide numerically exact answers?
More specifically, to what extent is the \emph{reasoning} of a floating-point LP-based branch-and-bound solver, which is composed of a large number of critical decisions during the search, correct in an exact sense?
While studies exist that compare the final result of numeric and exact solvers, e.g., in~\cite{CookKochSteffyWolter2013}, we are not aware of any analysis that looks deeper at the branch-and-bound process.
This is particularly interesting for the large number of cases where the final answer of a floating-point solver is correct: Is this due to the fact that all individual decisions were correct, or because all incorrect decisions did not affect or compromise the final result?




In this paper, we try to give empirical answers to these questions by performing a first such \emph{a posteriori} analysis of a pure LP-based branch-and-bound process.
To this end, we inspect all leaves of the branch-and-bound tree and check if the leaf is pruned correctly because
\begin{enumerate*}[label=(\alph*)]
    \item its LP relaxation is infeasible,
    \item its dual bound proves that the leaf does not contain improving solutions, or
    \item the LP solution is primal feasible
\end{enumerate*}.
To test the correctness of each decision, we apply a hierarchy of techniques from the literature such as safe bounding by directed rounding~\cite{NeumaierShcherbina04SafeBounds}, rational reconstruction by continued fraction approximations~\cite{UrsicPatarra83ExactSolutionOfSystems,Wan06AlgorithmTSolveLSexactly,Pan11NearlyOptimalSolutions,Saundersetall11IterativeRefinement}, exact LU factorization of simplex bases, and, if necessary, round-off-error-free LP solving~\cite{GleixnerSteffy2019,GleixnerSteffyWolter16IterativeRefinement}.
%
We base our experiments on the open-source MIP solver \scip~\cite{SCIP9} and use two curated test sets from the literature with and without numerical challenges~\cite{CookKochSteffyWolter2013}.
Our code base is publicly available.\footnote{\url{https://github.com/alexhoen/bnbanalyzer}}

The paper is organized as follows.
In \Cref{sec::description}, we briefly introduce the general concept of a \bandb solver, provide a classification of numerical errors, and survey the techniques to check the correctness of floating-point decisions. 
In \Cref{sec::experiments}, we present the results of our experiments and analyze the different errors encountered in the \fpsolver.
In \Cref{sec::conclusion}, we summarize our results and give an outlook on their implications for future research.

%% file: sections/basics.tex
\tikzset{
	tree node/.style = {align=center, inner sep=0pt, font = \scriptsize},
	S/.style = {draw, circle, minimum size = 8mm, top color=white, bottom color=white},
	tree node label/.style={font=\scriptsize},
}
\forestset{
	declare toks={left branch prefix}{},
	declare toks={right branch prefix}{},
	declare toks={left branch suffix}{},
	declare toks={right branch suffix}{},
	tree node left label/.style={
		label=170:#1,
	},
	tree node right label/.style={
		label=10:#1,
	},
	maths branch labels/.style={
		branch label/.style={
			if n=1{
				edge label={node [left, midway] {$\forestoption{left branch prefix}##1\forestoption{left branch suffix}$}},
			}{
				edge label={node [right, midway] {$\forestoption{right branch prefix}##1\forestoption{right branch suffix}$}},
			}
		},
	},
	text branch labels/.style={
		branch label/.style={
			if n=1{
				edge label={node [left, midway] {\foresteoption{left branch prefix}##1\forestoption{left branch suffix}}},
			}{
				edge label={node [right, midway] {\forestoption{right branch prefix}##1\forestoption{right branch suffix}}},
			}
		},
	},
	text branch labels,
	set branch labels/.style n args=4{%
		left branch prefix={#1},
		left branch suffix={#2},
		right branch prefix={#3},
		right branch suffix={#4},
	},
	set maths branch labels/.style n args=4{
		maths branch labels,
		set branch labels={#1}{#2}{#3}{#4},
	},
	set text branch labels/.style n args=4{
		text branch labels,
		set branch labels={#1}{#2}{#3}{#4},
	},
	branch and bound/.style={
		/tikz/every label/.append style=tree node label,
		maths branch labels,
		for tree={
			tree node,
			S,
			math content,
			s sep'+=20mm,
			l sep'+=5mm,
			thick,
			edge+={thick},
		},
		before typesetting nodes={
			for tree={
				split option={content}{:}{tree node left label,content,tree node right label,branch label},
			},
		},
		where n children=0{
			tikz+={
				\draw [thick]  ([yshift=-10pt, xshift=-2.5pt].south west) -- ([yshift=-10pt, xshift=2.5pt].south east);
			}
		}{},
	},
}

\section{Numerical Errors in LP-based Branch and Bound}
\label{sec::description}

A \emph{mixed integer program} (MIP) is given in the form
\begin{align*}
  \min\; & c^{T}x\\
  \text{s.t. }\; & Ax \geq b,\\
  & \ell \leq x \leq u,\\
  & x_i \in \mathbb{Z}\;\; \text{for all } i \in \mathcal{I},
\end{align*}
where we assume rational input data $A\in \mathbb{Q}^{m \times n}$, $c\in\mathbb{Q}^{n}$, $u,\ell\in (\mathbb{Q}\cup\{\pm\infty\})^{n}$, and $b \in \mathbb{Q}^m $.
The index set $\mathcal{I} \subseteq \{1,\dots,n\}$ defines which decision variables are required to be integer.

MIPs are typically solved by variants of \emph{LP-based branch and bound}, which dates back to \cite{LandDoig60BranchAndBound} and still forms the backbone of today's most competitive MIP solvers.
In \Cref{sec::bas::branchbound}, we describe briefly the procedure of an LP-based \bandb algorithm.
For a more detailed description we refer to \cite{Achterberg07Thesis}.
In \Cref{sec::bas::resolving}, we give an overview of the techniques we use to obtain an exact evaluation of a node that was solved in floating-point arithmetic.
In \Cref{sec::bas::error}, we categorize how numerics can impact the critical \bandb decisions in a floating-point solver and how incorrect decisions affect the overall result.

\subsection{Basics of the LP-based Branch-and-Bound Algorithm}
\label{sec::bas::branchbound}

The first step in an LP-based \bandb algorithm is to \emph{relax} the problem by dropping all integrality constraints.
The resulting LP relaxation can be solved by an LP solver. 
If the LP solution is \emph{integer feasible}, i.e., all integer variables have an integral solution value, the problem is solved.
If the LP relaxation is detected to be infeasible, then also the original MIP is infeasible.
Otherwise, the algorithm \emph{branches}:
It chooses an integer variable $x_i$ with a fractional value $\hat x_i$ in the solution of the LP relaxed and creates two new sub-problems
by taking the original problem and adding the constraint $x_i \geq \lceil \hat x_i \rceil$ and $x_i \leq \lfloor \hat x_i \rfloor$, respectively.
This ensures that the optimal solution is preserved in at least one of the two sub-problems and the current LP solution is no longer feasible in both sub-problems.
This process of LP solving and branching is iterated recursively.

If during this process an integer feasible solution is found its objective value can be used to update the \emph{primal bound}, which is the objective value of the best known solution encountered so far.
This primal bound can be used to \emph{prune}, i.e., remove from further consideration all open sub-problems for which the\emph{dual bound}, i.e., the objective value of the sub-problem's LP relaxation, is greater than or equal to the primal bound.

\begin{figure}[ht]
    \centering

    \begin{forest}
        branch and bound,
        where level=1{
            set branch labels={x_2\leq}{}{x_2\geq}{},
        }{
            if level=2{
                set branch labels={x_3\leq}{}{x_3\geq}{},
            }{},
        }
        [
        $1.5$:root:\{0;$\frac{1}{2}$;1\}
            [
            $2.5$:N_1:\{1;0;$\frac{1}{2}$\}:0]
            [
            $1.5$:N_2:\{0;1;$\frac{1}{2}$\}:1
                [
                $\infty$:N_3:LP infeas.:0]
                [2:N_4:\{0;1;1\}:1]
            ]
        ]
    \end{forest}
    \caption{An example for a \bandb tree on a minimization problem.}
    \label{fig::tree}
\end{figure}
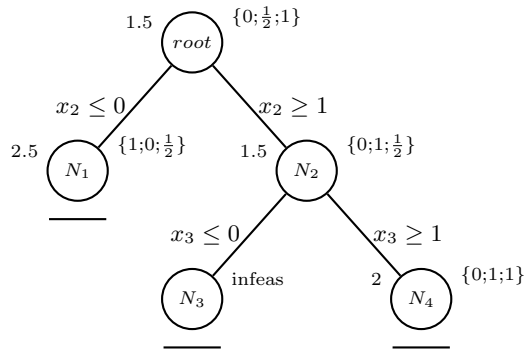

This recursive procedure produces a search tree with sub-problems as nodes.
A node is called \emph{leaf} if it is not further branched on because it is     
\begin{enumerate*}[label=(\alph*)]
    \item LP infeasible, 
    \item it was pruned, or
    \item the LP solution is integer feasible.
\end{enumerate*}
\Cref{fig::tree} presents an example of a \bandb tree for a minimization problem.
The solution of the relaxed LP at the node is displayed slightly above and to the right of the node, while its objective value is positioned on the left side.
The leaf $N_4$ is integer feasible and appears to be also optimal.
The leaf $N_3$ is LP infeasible and the node $N_1$ is cut off since branching on $N_1$ does not improve the solution of $N_4$ anymore.

Note that although the heuristic decision on which fractional variable to branch may be affected by numerical calculations, any sequence of branching decisions yields a correct branch-and-bound process.
The same holds for the heuristic decision in which order to process open sub-problems: any order yields a correct \bandb algorithm.
By contrast, the decision of whether to remove a sub-problem from consideration is both affected by numerics and can critically affect the correctness of the final result.
In the next section, we focus on these critical decisions.

\subsection{Types of Numerical Errors in LP-based Branch and Bound}
\label{sec::bas::error}

In our analysis, we distinguish the following types of incorrect decision that can be taken by a floating-point solver due to numerical errors.
All of these errors occur at a leaf of the \bandb tree.
%

\paragraph{Solution errors.}
    By design, a \fpsolver may produce approximate solutions with slight violations of integrality and constraint feasibility. A simple post-processing step is to first round the value of all integer variables to the nearest integer and second, if continuous variables are present, to compute rational values for these such that all constraints are satisfied exactly.
    If this is certifiably not possible, i.e., if the \fpsolver accepted an assignment for the integer variables as feasible that cannot be completed to an exactly feasible solution, we refer to this as a \emph{solution error}.
    
    We distinguish two types of solution errors.
    If the LP relaxation at this node is infeasible in exact arithmetic, then no solution was discarded by removing the current leaf.
    We refer to this error as a \emph{weak solution error}.
    Otherwise, the primal solution of the LP relaxation does not satisfy integer feasibility, but the sub-problem below the current leaf node may still contain improving solutions.
    In this case, only a lower bound can be derived and we refer to this situation as a \emph{strong solution error}.

\paragraph{Bound errors.}
    If a \fpsolver uses a lower bound to prune a node, this decision can be verified \emph{a posteriori} by computing solving the LP relaxation exactly or computing a safe dual bound that is valid in exact arithmetic, by techniques described in \Cref{sec::bas::resolving}.
    If such a safe bound is greater than or equal to the current primal bound, then the decision to prune the node was correct.
    Otherwise, if the exact objective value of the LP relaxation is below the primal bound, we label this decision as a \emph{bound error}.
    Again, we distinguish two cases.  
    If later during the solving process an exactly feasible solution is found that improves the primal bound and in hindsight justifies pruning the leaf, we call the error a \emph{weak bound error}, otherwise a \emph{strong bound error}.

\paragraph{Gap errors.}
    For the third type of error, consider the situation that a \fpsolver declares a node to be integer feasible and produces an approximate solution that can be converted to an exactly feasible solution after rounding the integer assignment.
    Although we record no solution error here, the decision to terminate the search below this node may still be incorrect if the objective value of the converted, exact solution does not match the dual bound derived by solving the LP relaxation exactly.
    We refer to this error as a \emph{gap error}, and as above we distinguish between \emph{weak gap errors} that are justified in hindsight and \emph{strong gap errors}. 

\paragraph{Infeasibility errors.}
    Finally, if the \fpsolver declares a leaf infeasible but there exists a solution to the LP relaxation that is exactly feasible we refer to this situation as an \emph{infeasibility error}.


~\\

\noindent
Even when some of the errors listed above occur, the overall result of the \fpsolver may still be correct.
Specifically, weak bound and weak gap errors do not compromise the correctness of the final result since the corresponding decisions are justified in hindsight.
We still record these errors, because we are interested not only in the correctness of the final result but in the correctness of the \fpsolver's reasoning that led to the result.

Also in the presence of strong bound or strong gap errors, the solution or at least its objective value~$z^*$ returned by the \fpsolver may still be optimal.
Because our analysis produces a safe dual bound at all occurrences of errors outlined above, we can then take the minimum of all these dual bounds and derive a global lower bound~$\hat z < z^*$ for the true optimal objective value.
This certifies that the true optimal objective value must lie in the interval~$[\hat z,z^*]$.
If infeasibility errors occur, i.e., if a leaf is removed although the LP relaxation is feasible in exact arithmetic, we can use the exact objective value of the LP relaxation as a valid dual bound.

Finally, note that in the presence of solution errors or gap errors, the \fpsolver will work with an incorrect, usually too optimistic primal bound during part of the search.
When checking for bound errors, however, we always take the perspective of the \fpsolver, i.e., we treat the primal bound registered in the \fpsolver at that time as valid, and check whether the pruning decision is justified with respect to this primal bound.
This helps us to separate more clearly and quantify more precisely the different types of incorrect reasoning in our computational analysis later.

\subsection{Checking Numerical Correctness of Floating-point Decisions}
\label{sec::bas::resolving}

A straightforward way to evaluate whether and which of the errors described above occurs at a leaf node is to solve the LP relaxation once again with an exact LP solver and compare the results.
For most but small instances, this approach is prohibitively slow.
Hence, we use a hierarchy of methods that may be able to verify or refute the correctness of leaf decisions faster, described in the following.

\paragraph{Safe bounding.}
Consider the dual of the LP relaxation,
\[
\max \big\{ b^T y + \ell^T r^+ - u^T r^- :
A^T y + r^+ - r^- = c,\,
y, r^+, r^- \geq 0 \big\},
\]
and assume we have an approximate dual solution $(\hat{y}, \hat{r}^+, \hat{r}^-)$ with dual residual error 
\begin{align}
    \hat\varepsilon := c - A^T\hat{y} - \hat{r}^+ + \hat{r}^- \in \mathbb{Q}^n.
\end{align}
If we compute, in exact arithmetic,
\begin{align}
    r^+_i := \hat{r}^+_i + \max\{\hat\varepsilon_i,0\}
    \textnormal{ and } 
    r^-_i := \hat{r}^-_i - \min\{\hat\varepsilon_i,0\} 
\end{align} 
for all variables $i \in \{1,\ldots,n\}$, then $\{\hat{y},r^+,r^-\}$ is an exactly feasible dual solution with objective value $b^T\hat{y} + \ell^T r^+ - u^T r^-$.
This modification thus ensures dual feasibility at the cost of reducing the objective value by $\ell^T(\hat{r}^+- r^+)+ u^T (r^- - \hat{r}^-)$. 
The quality of the resulting safe dual bound is impacted by the variable bounds $u,\ell$, which, ideally, should be as tight as possible to minimize the amount by which the approximate dual bound is weakened.
If infinite variable bounds are used during the computation, the resulting safe bound becomes negative infinity.
This technique is first explained in~\cite{NeumaierShcherbina04SafeBounds} in a version using interval arithmetic and discussed in more detail in \cite{SteffyWolter2011}.
Safe bounding can be applied analogously to validate a Farkas proof for infeasible LPs.

\paragraph{Rational reconstruction.}
In the hope of recovering an exact rational solution value from an approximate floating-point value, we use continued fractions approximation, which can be computed by the extended Euclidean algorithm.
This way, we can round floating-point values to a nearby rational number with a limited denominator.
This method is applied to both primal and dual LP solutions, as well as to Farkas proofs.
Originally, this technique is described and improved in \cite{Pan11NearlyOptimalSolutions,Saundersetall11IterativeRefinement,UrsicPatarra83ExactSolutionOfSystems,Wan06AlgorithmTSolveLSexactly}.
We rely on an implementation in the exact LP solver \soplex~\cite{GleixnerSteffy2019,GleixnerSteffyWolter16IterativeRefinement}.

\paragraph{Factorization.}
If a simplex basis is available for the LP relaxation solution, we perform an exact rational LU factorization in order to compute the exact rational solution associated to the basis reached by the floating-point LP solver. 
Again, we use the implementation in \soplex~\cite{GleixnerSteffy2019,GleixnerSteffyWolter16IterativeRefinement}.
This step is only skipped for LP infeasible nodes since no basis can be obtained from the LP solver.

\paragraph{Exact LP solving.}

The last and most expensive option to derive the most accurate lower bound on the objective is to call a round-off-error-free LP solver. 
The open-source solver \soplex provides a configuration to perform such an exact rational solve~\cite{GleixnerSteffy2019,GleixnerSteffyWolter16IterativeRefinement}.
We refer to this fallback as \exactsoplex.
It is performed only if all techniques described above are not applicable or fail to decide the status of a leaf.

Exact LP solving is also used to test whether an approximate primal solution can be converted to an exact primal solution by fixing the integer variables to the rounded values of the floating-point solution and solving the LP over the continuous variables in rational arithmetic with \exactsoplex.
%
For each primal solution value, it is ensured to lie within its lower and upper bounds. 

%


%% file: sections/experiments.tex
\section{Computational Study}
\label{sec::experiments}

The guiding questions for the design of our experiments are the following:
Assuming the final result produced by a floating-point \bandb solver is correct, can we verify that the also reasoning of the individual decisions of the solver is correct?
If not, how many of the critical decisions would require re-evaluation, and which of the errors defined in \Cref{sec::bas::error} occur at which frequency?
Furthermore, how effective are the different post-processing techniques outlined in \Cref{sec::bas::resolving}.
Before evaluating our experiments in \Cref{sec::exp::exp}, we first describe the setup and specify the software used in \Cref{sec::exp::setup}.

\subsection{Experimental Setup}
\label{sec::exp::setup}

As \bandb framework, we use the open-source MIP solver \scipv\cite{Achterberg07Thesis,SCIP9} \sciphash\xspace 
and configure \scip via parameters such that it mimics the behavior of the pure LP-based \bandb
algorithm described in~\cite{Cooketal11ExactRationalMIPSolver}.
We use the open-source solver \soplexv \soplexhash\xspace both as the underlying floating-point LP solver in \scip, and as an exact LP solver during our \emph{a posteriori} analysis.

\paragraph{Preprocessing.}
Before passing the instances to \scip and starting the actual \bandb process, we apply simple presolving steps~\cite{presolving_achterberg,BrearleyMitra75AnalysisOfMPP}. 
First, we perform some \emph{model cleanup}, i.e., we delete redundant constraints and convert singleton rows to variable bounds. 
Further, we apply a limited form of \emph{constraint propagation} in rational arithmetic in order to reduce the number of infinite bounds and help increase the success rate of safe dual bounding.

\paragraph{Event system.}
%
In order to track \scip's \bandb process, we register so-called events that inform us every time a leaf is reached.
In detail, our implementation uses the following events.
        Every time  \scip finds a solution, it triggers the event \bestsol. 
        A \nodeinfeas event is triggered both when the LP relaxation at the current leaf is deemed infeasible, or when a node is pruned.
        We access the LP solving status to differentiate whether the relaxed LP is infeasible or feasible and eventually pruned.
        Note that when \scip creates new nodes it adds them to a queue and sets the local lower bound to the objective value of the LP relaxation of its parent.
        After a new solution is found, \scip immediately removes all nodes from the queue whose lower bound is worse than the objective of the newly found solution. 
        Each of these removed leaves triggers a \nodedel event.
        Finally, every integer-feasible node triggers a \nodefeas event.

All node events come with a list of branching decisions.
Additionally for the  \nodefeas and \nodeinfeas events, also the primal and dual LP solution, respectively the Farkas proof \cite{Farkas1902}, and the simplex basis are available.
When the \nodedel event is triggered, we recompute the dual LP solution by solving the LP again in floating-point arithmetic with \soplex. 
At the \bestsol event, \scip only provides the primal solution.

\paragraph{Test sets.}
We consider two test sets curated in \cite{Cooketal11ExactRationalMIPSolver} for benchmarking an exact MIP solver.
The \fpeasy test set consists of 57 instances that were found easy by the floating-point solver \scip.
The \numdiff test set consists of 50 instances and was compiled to be numerically challenging due to large coefficient ranges or poor conditioning. 
For a detailed description we refer to \cite{Cooketal11ExactRationalMIPSolver}, which also lists on which instances the results of floating-point \scip deviate from the exact rational results.

To provide more stable results, we perform our experiments on each instance twice: one time on the original instance and another time after permuting the order of variables and constraints.
For the \numdiff test set instance \solver{normalized-aim-200-1\_6-ye} we use \solver{aim-200} as an abbreviation.

The code base used for our experiments is publicly available at \url{https://github.com/alexhoen/bnbanalyzer}.
%
All experiments were carried out on identical machines with Intel(R) Xeon(R) CPU E7-8880 v4 @ 2.20 GHz with a time limit of 10800 seconds for each instance, which includes both the solving time of the \fpsolver \scip as well as the evaluation of the leaves.


\subsection{Evaluation of Branch-and-Bound Decisions in \scip}
\label{sec::exp::exp}

\Cref{tab::stats} provides a first overview of the results for the \fpsolver \scip on both test sets \fpeasy and \numdiff, summarizing how many instances were correctly solved and on how many \scip returns an inexact solution or even an incorrect status. 
The columns ``correct'' state the number of instances where no errors occurred as defined in \Cref{sec::bas::error}:
All solutions are either exactly integer feasible or convertible to such, all infeasible leaves are infeasible in rational arithmetic, and any node pruning is also justified in exact arithmetic.
By contrast, the columns ``fails'' state the number of instances where at least one such error occurred and solving finished within the time limit. 
The column ``primal'' reports the number of instances with a solution error (see \Cref{sec::bas::error}), the column ``dual'' reports the number of instances with a bound, gap, or infeasibility error.
%
Note that for the two \solver{MIPLIB}-instances (\solver{30:70:4\_5:0\_95:100} and \solver{npmv07}) \scip timed out before reaching a leaf or producing a solution and therefore did not issue any event.
Since no wrong decision was made these instances are labeled as correct. 


\begin{table}
    \begin{tabular*}{\textwidth}{@{}l@{\;\;\extracolsep{\fill}}cccc ccc cccc  }
        \toprule
        \multicolumn{3}{c}{} & \multicolumn{4}{c}{within time limit} & \multicolumn{2}{c}{time out} \\
        \cmidrule(lr){4-7} \cmidrule(lr){8-9}
        \multicolumn{2}{l}{test set} & instances & correct & fails & primal & dual & correct & fails\\
        \midrule
        \fpeasy & original & 57 & 49 & 2 & 0 & 2 & 6 & 0\\
                & permuted & 57 & 49 & 2 & 0 & 2 & 6 & 0\\
        \midrule
        \numdiff & original & 50 & 11 & 20 & 12 & 9 & 14 & 5\\
                 & permuted & 50 & 13 & 18 & 12 &  9 & 13 & 7\\
        \bottomrule
    \end{tabular*}
    \caption{Aggregate statistics on instances with/without incorrect decisions in floating-point \scip.}
    \label{tab::stats}
\end{table}

Regardless of the permutation, for 55 of the 57 instances of \fpeasy, we could verify that \scip followed a fully correct solving process.
On 49~instances it terminated with an exact optimal solution, while on 6 instances \scip hit the time limit.
Only for two instances, we encountered dual fails.
On one instance (\solver{vpm2}) the result is still correct, on the other instance (\solver{dano3\_4}), the result is slightly suboptimal.
This can also be seen in \Cref{tab::diff_scip_exactscip}, where we compare the results of floating-point \scip to the results of a numerically exact version of \scip~\cite{Cooketal11ExactRationalMIPSolver,EiflerGleixner2022_Acomputational,eiflergleixner2023safeverifiedgomorymixed}, on the instances where \scip made at least one wrong decisions.

As expected, verifying the results of the \numdiff test set exhibits more errors. 
Only on 25 (original) and 27 (permuted) of the 50 instances, respectively, \scip followed a fully correct solution process.
On each seed, 14 of these instances hit the time limit, so only for 11 respectively 13 instances \scip terminated with a verifiably exact optimal solution.
Still, on both test sets, a notable amount of the floating-point solves exhibits no errors.

For a more detailed analysis, \Cref{tab::error::fpeasy,tab::error::numdiff,tab::error::numdifftlim} include all instances labeled as ``fail'' in \Cref{tab::stats} and those instances for which \scip produced a ``correct'' result but made incorrect decisions during the solving process.
All instances where \scip made a wrong decision on \fpeasy are listed in \Cref{tab::error::fpeasy}.
Instances of \numdiff that encountered errors and are solved within the time limit are listed in \Cref{tab::error::numdiff}.
A complete overview \numdiff instances with wrong decisions where \scip timed out are listed in \Cref{tab::error::numdifftlim}.
The column ``leaves'' lists the total number of processed leaves of the tree; the remaining columns report the number of leaves with errors according to the definition in \Cref{sec::bas::error}, where ``W'' and ``S'' stands for weak and strong errors, respectively.

\begin{table}
    \begin{tabular*}{\textwidth}{@{}l@{\;\;\extracolsep{\fill}}rrrrrrrrr}
        \toprule
        \multicolumn{3}{c}{} & \multicolumn{2}{c}{Sol} & \multicolumn{2}{c}{Bound} & \multicolumn{2}{c}{Gap} & Inf \\
        \cmidrule(lr){4-5} \cmidrule(lr){6-7} \cmidrule(lr){8-9}
        \multicolumn{1}{c}{Instance} & Perm & leaves & W & S & W & S & W  & S &  \\
        \midrule
        \multirow{2}{*}{\solver{dano3-4}} 
            & 0 & 29 & 0 & 0 & 0 & 0 & 0 & 1 & 0\\
            & 1 & 29 & 0 & 0 & 0 & 0 & 0 & 1 & 0\\
        \midrule
        \multirow{2}{*}{\solver{vpm2}}    
            & 0 & 688\,018 & 0 & 0 & 5  & 34  & 0 & 0 & 0 \\
            & 1 & 649\,085 & 0 & 0 & 11 & 167 & 0 & 0 & 0\\
        \bottomrule
    \end{tabular*}
    \caption{Analysis of number of leaves with incorrect decisions on the \fpeasy instances within the time limit.}
    \label{tab::error::fpeasy}
\end{table}

Notably, only the instances \solver{ns1859355} and \solver{vpm2} show a correct result despite wrong decisions during the solving process.
On \fpeasy the vast majority of instances in \Cref{tab::error::fpeasy} are solved without wrong decisions.
SCIP incorrectly evaluates nodes on only two instances, affecting one node in \solver{dano3-4} and 39 of 178 nodes in \solver{vpm2}, respectively.
Though 39 or 178 wrong decisions on \solver{vpm2} seems to be a notable amount of wrong decisions, compared to the overall evaluated nodes on the entire \fpeasy test set the wrong decisions on these two instances make up a tiny portion.
\begin{table}
    \begin{tabular*}{\textwidth}{@{}l@{\;\;\extracolsep{\fill}}rrrrrrrrr}
        \toprule
        \multicolumn{3}{c}{} & \multicolumn{2}{c}{Sol} & \multicolumn{2}{c}{Bound} & \multicolumn{2}{c}{Gap} & Inf \\
        \cmidrule(lr){4-5} \cmidrule(lr){6-7} \cmidrule(lr){8-9}
        \multicolumn{1}{c}{Instance} & Perm & leaves & W & S & W & S & W  & S  &  \\
        \midrule
        \solver{alu10\_1} & 0 & 883 & 2 & 0 & 0 & 0 & 0 & 0 & 0 \\
        \midrule
        \multirow{2}{*}{\solver{alu10\_7}}
            & 0 & 811 & 0 & 1 & 0 & 0 & 0 & 0 & 0 \\
            & 1 & 734 & 0 & 3 & 0 & 0 & 0 & 0 & 0 \\
        \midrule
        \multirow{2}{*}{\solver{alu10\_8}}
            & 0 & 4\,254 & 0 & 3 & 0 & 0 & 0 & 0 & 0 \\
            & 1 & 4\,486 & 1 & 3 & 0 & 0 & 0 & 0 & 0 \\
        \midrule
        \multirow{2}{*}{\solver{alu10\_9}}
            & 0 & 6\,665 & 1 & 3 & 0 & 0 & 0 & 0 & 0 \\
            & 1 & 7\,675 & 1 & 5 & 0 & 0 & 0 & 0 & 0 \\
        \midrule
        \solver{alu16\_1} 
            & 0 & 958\,643 & 0 & 0 & 0 & 0 & 0 & 0 & 46 \\
        \midrule
        \multirow{2}{*}{\solver{alu16\_7}}
            & 0 & 694 & 0 & 3 & 0 & 0 & 0 & 0 & 0 \\
            & 1 & 1\,307 & 0 & 4 & 1 & 1 & 0 & 0 & 0 \\
        \midrule
        \multirow{2}{*}{\solver{alu16\_8}}
            & 0 & TL & 0 & 8 & 0 & 20 & 0 & 0 & 3 \\
            & 1 & 219\,362 & 0 & 3 & 0 & 2 & 0 & 0 & 0 \\
        \midrule
        \multirow{2}{*}{\solver{alu16\_9}}
            & 0 & 199\,225 & 0 & 4 & 0 & 0 & 0 & 0 & 0 \\
            & 1 & 44\,961 & 0 & 13 & 0 & 0 & 0 & 0 & 0 \\
        \midrule
        \multirow{2}{*}{\bernd}
            & 0 & 31\,081 & 0 & 4 & 0 & 0 & 0 & 0 & 0 \\
            & 1 & 23\,144 & 0 & 5 & 0 & 0 & 0 & 0 & 0 \\
        \midrule
        \multirow{2}{*}{\solver{dfn6\_load}}
            & 0 & 11\,210 & 2 & 6 & 9 & 0 & 0 & 0 & 0 \\            
            & 1 & 3\,064 & 0 & 6 & 9 & 0 & 0 & 0 & 0 \\
        \midrule
        \multirow{2}{*}{\solver{aim-200}}
            & 0 & 7\,414 & 0 & 0 & 0 & 0 & 0 & 0 & 1 \\
            & 1 & 9\,534 & 0 & 0 & 0 & 0 & 0 & 0 & 1 \\
        \midrule
        \multirow{2}{*}{\solver{neos-1053591}}
            & 0 & 225\,381 & 0 & 2 & 0 & 0 & 0 & 0 & 0 \\
            & 1 & TL & 1 & 1 & 0 & 0 & 0 & 0 & 0 \\
        \midrule
        \multirow{2}{*}{\solver{neos-1062641}}
            & 0 & 92 & 0 & 1 & 0 & 0 & 0 & 0 & 0 \\
            & 1 & 46 & 0 & 1 & 0 & 0 & 0 & 0 & 0 \\
        \midrule
        \solver{neos-1603965} & 0/1 & 1 & 0 & 1 & 0 & 0 & 0 & 0 & 0 \\
        \midrule
        \multirow{2}{*}{\solver{ns1629327}}
            & 0 & 9\,985 & 0 & 0 & 0 & 0 & 0 & 2 & 0 \\
            & 1 & 11\,788 & 0 & 0 & 0 & 2 & 0 & 2 & 0 \\
        \midrule
        \multirow{2}{*}{\solver{ns1859355}}
            & 0 & 12\,581 & 0 & 0 & 0 & 0 & 1 & 0 & 0 \\
            & 1 & 6\,434 & 0 & 0 & 0 & 0 & 2 & 0 & 0 \\
        \midrule
        \solver{ns1866531} & 0/1 & 1 & 0 & 1 & 0 & 0 & 0 & 0 & 0 \\
        \midrule
        \multirow{2}{*}{\solver{prodplan2}}
            & 0 & 3 & 0 & 0 & 0 & 1 & 1 & 1 & 0 \\
            & 1 & 27 & 0 & 0 & 3 & 0 & 2 & 1 & 0 \\
        \midrule
        \multirow{2}{*}{\solver{tkat3K}}
            & 0 & 4\,367 & 0 & 0 & 11 & 1 & 0 & 0 & 0 \\
            & 1 & 8\,335 & 0 & 0 & 5 & 2 & 0 & 0 & 0 \\
        \midrule
        \multirow{2}{*}{\solver{tkat3T}}
            & 0 & 14\,604 & 0 & 0 & 27 & 1 & 0 & 0 & 0 \\
            & 1 & 5\,219 & 0 & 0 & 4 & 4 & 0 & 0 & 0 \\
        \midrule
        \multirow{2}{*}{\solver{tkat3TV}}
            & 0 & 4\,546 & 0 & 0 & 10 & 7 & 0 & 0 & 0 \\
            & 1 & 19\,032 & 0 & 0 & 28 & 5 & 0 & 0 & 0 \\
        \midrule
        \multirow{2}{*}{\solver{tkatTV5}}
            & 0 & 8\,459 & 0 & 0 & 15 & 26 & 0 & 0 & 0 \\
            & 1 & 18\,999 & 0 & 0 & 1 & 17 & 0 & 0 & 0 \\
        \bottomrule
    \end{tabular*}
    \caption{Analysis of number of leaves with incorrect decisions on the \numdiff instances within the time limit.}
    \label{tab::error::numdiff}
\end{table}

\begin{table}
    \begin{tabular*}{\textwidth}{@{}l@{\;\;\extracolsep{\fill}}rrrrrrrrr}
        \toprule
        \multicolumn{3}{c}{} & \multicolumn{2}{c}{Sol} & \multicolumn{2}{c}{Bound} & \multicolumn{2}{c}{Gap} & Inf \\
        \cmidrule(lr){4-5} \cmidrule(lr){6-7} \cmidrule(lr){8-9}
        \multicolumn{1}{c}{Instance} & Perm & leaves & W & S & W & S & W  & S  & \\
        \midrule
        \multirow{2}{*}{\solver{alu16\_5}} 
            & 0 & 1\,342\,439 & 0 & 0 & 0 & 0 & 0 & 0 & 11 \\
            & 1 & 280\,574 & 0 & 0 & 0 & 0 & 0 & 0 & 875 \\
        \midrule
        \solver{dfn6fp\_load} & 1 & 585 & 0 & 2 & 3 & 4 & 0 & 0 & 0 \\
        \midrule
        \multirow{2}{*}{\solver{ns2080781}} 
            & 0 & 701\,865 & 7 & 3 & 156 & 72 & 0 & 0 & 0 \\
            & 1 & 52\,258 & 0 & 4 & 0 & 0 & 0 & 0 & 0 \\
        \midrule
        \solver{ns1925218} & 1 & 22\,556 & 0 & 0 & 0 & 0 & 0 & 0 & 1 \\
        \midrule
        \solver{prodplan1} & 0 & 36 & 0 & 0 & 0 & 1 & 0 & 0 & 0 \\
        \midrule
        \solver{ran14x18-}
            & 0 & 1\,268\,265 & 0 & 12 & 1 & 0 & 0 & 0 & 0 \\
             \solver{-disj-8}
            & 1 & 1\,288\,277 & 0 & 6 & 2 & 2 & 0 & 0 & 0 \\
        \midrule
        \solver{neos-799716} & 1 & 1\,095 & 0 & 0 & 0 & 0 & 0 & 0 & 1 \\
        \bottomrule
    \end{tabular*}
    \caption{Analysis of number of leaves with incorrect decisions on the \numdiff instances with time out.}
    \label{tab::error::numdifftlim}
\end{table}

\begin{table}
    \begin{tabular*}{\textwidth}{@{}l@{\;\;\extracolsep{\fill}}| rr | l}
        \toprule
         \multicolumn{1}{c}{}& \multicolumn{2}{c}{\scip}  & \multicolumn{1}{c}{range of exact solution}\\
        \cmidrule{2-3}
        \multicolumn{1}{c}{Instance} & Seed 0 & \multicolumn{1}{c}{Seed 1} &  \\
        \midrule
        \solver{dano3\_4} & 576.435224722083 & 576.435224722083 & 576.4352247072 \\
        \solver{vpm2} & 13.75 & 13.75 & 13.75 \\
        \midrule
        \solver{alu10\_1} & 86 & inf & inf  \\
        \solver{alu10\_7} & 83 & 83 & [1243230.17385925,$\infty$)  \\
        \solver{alu10\_8} & 84 & 84 & [169.4062227788,$\infty$)  \\
        \solver{alu10\_9} & 84 & 83 & [1166.7499882579,3726773] \\
        \solver{alu16\_1} & 84 & inf  & inf \\
        \solver{alu16\_7} & 79 & 79 & [315.4000244141,$\infty$)  \\
        \solver{alu16\_8} & x & 79 & [2256.0481863315,$\infty$) \\
        \solver{alu16\_9} & 79 & 79 & [3801.0189752579,$\infty$)  \\
        \solver{bernd2}   & 11209.0573895658 & 11209.0573899187 & 113091.469015879 \\
        \solver{dfn6\_load} & 3.743842258432 & 3.743842258432 & [3.7683622392,4.4007262645] \\
        \solver{aim-200} 
                          & inf & inf & 200 \\
        \solver{neos-1053591} & -3662.9144 & time limit & -3662.9144\\
        \solver{neos-1062641} & 0 & 0 & 0  \\
        \solver{neos-1603965} & 619244367.662956 & 619244367.662956 & [619246130.539662,$\infty$) \\
        \solver{ns1629327}  & -10.9803191329325 & -10.9803191329325 & -10.9803191329  \\
        \solver{ns1866531}  & 0 & 0 & 10 \\
        \solver{prodplan2}  & -239399.435140992 & -239399.4351411 & -239399.435141407 \\
        \solver{tkat3K}     & 4772818.1 & 4772818.1 & 4772818.1 \\
        \solver{tkat3T}     &  5564891.75 & 5564891.75 & 5564891.75  \\
        \solver{tkat3TV}    & 8388398.65 & 8388398.65 & 8388398.65  \\
        \solver{tkatTV5}    & 28117644.225 &28117644.225 	 & 28117644.225\\
        \bottomrule
    \end{tabular*}
    \caption{Comparison of the results of floating-point \scip and rational solutions. The exact results were generated with \exactscip \exactsciphash\xspace and a time limit of 18800 seconds.}
    \label{tab::diff_scip_exactscip}
\end{table}

As expected, on the numerically more challenging \numdiff test set, more wrong decisions are made.
First, let us focus on the instances with dual fails.
Besides the instances \solver{alu16\_8}, \solver{tkatTV5}, and  \solver{tkat3TV}, all instances show at most two leaves with strong bound or gap errors.
On the instances, \solver{tkat3*}, \solver{ns1629327} and \solver{ns1859355}, more weak and strong bound or gap errors appear, and \solver{alu10\_9} produces 46 infeasibility errors.
Nevertheless, these all make up only a small fraction of the total number of evaluated nodes.

On the \solver{alu*}-instances (except for \solver{alu16\_1}), as well as on the instances \dfnload, \solver{neos-1053591}, \solver{neos-1062641}, \solver{neos-1603965}, \solver{ns1866531} the solutions generated by \scip are slightly infeasible and can not be converted to exact solutions.
Especially the \solver{alu} instances are numerical difficult:
Either the instance is infeasible or there is a huge gap between the floating-point solution and the solution of \exactscip as can be seen in \Cref{tab::diff_scip_exactscip}.

One natural attempt to reduce such solution errors is to call \scip with a tightened primal feasibility tolerance. 
Per default, the feasibility tolerance is $10^{-6}$. 
We repeated our experiments but tightened the tolerance parameter ``numerics/feastol'' to $10^{-9}$.
The results show a significant reduction in solution errors.
On the instances \solver{alu10\_1}, \solver{alu10\_8}, \solver{alu10\_9}, \bernd, \solver{aim-200} and \solver{dfn6\_load}, \scip finds an exactly integer-feasible solution, or they can be converted to such, or \scip times out, but on the instances \solver{alu\_16\_7} or \solver{ns1859355}, the \fpsolver still produces slightly infeasible solutions.
On the \fpeasy instance \solver{dano3\_4} also all gap errors were removed, successfully closing the dual gap on this instance.


However, we also observe that the number of explored nodes is significantly increased.
This does not only lead to increased solving time, but can also create more numerically challenging leaves, hence potentially increasing the absolute number of incorrect decisions.
An example of this behavior is the instance \solver{blend2}, for which the result of \scip can no longer be verified after tightening the tolerance.


\subsection{Analysis of Different Post-processing Techniques}

Finally, let us have a look at the effectiveness and success rate of the different post-processing techniques described in \Cref{sec::bas::resolving}.
\Cref{tab::distribution} reports aggregate results over all leaves of a complete run over each test set (column ``leaves'') and categorizes them based on the techniques that were successful for verification.
The column ``floating-point'' lists the percentage of leaves that are verified automatically or through the use of safe bounding. 
The columns ``Reconstruct'', ``Factorize'', and ``\solver{Exact}'' likewise represent the percentage for the corresponding methods. 
The row ``unb'' reports the numbers for all instances that contain variables with unbounded domains.



\begin{table}\small
    \begin{tabular*}{\textwidth}{@{}l@{\;\;\extracolsep{\fill}}l|r  r  r  r  r r }
        \toprule
         \multicolumn{2}{c}{} & leaves & floating-point & Reconstruct & Factorize & {\solver{exact}} & error \\
        \midrule
        \multirow[t]{ 2}{*}{\fpeasy} 
        & all & 19952937 & 93.76\% & 2.93\% & 0.46\% & 2.85\% & 0.00020\%\\
        & unb & 424881 & 65.18\% & 0.15\% & 4.78\% & 29.89\% & 0.00024\% \\
       \multirow[t]{ 2}{*}{\numdiff} 
            & all & 8389838 & 99.69\%  & 0.0\%  & 0.04\% & 0.27\%  & 0.00543515\% \\
            & unb & 71221 & 94.89\%  & 0.02\%  & 1.22\% & 3.86\%  & 0.00140408\% \\
        \bottomrule
    \end{tabular*}
    \caption{Distribution on the different verification methods over all leaves on seed 0.}
    \label{tab::distribution}
\end{table}

Surprisingly, 93.76\% (\fpeasy) to 99.69\% (\numdiff) of the leaf decisions can be verified using floating-point arithmetic. 
Hence, the most successfull methods are also the ones that can be implemented most efficiently.
In relative terms, the percentage of erroneous leaf decisions is small on all test sets.

Safe bounding requires bounded variables and can therefore more reliably be applied to problems with bounded variables. 
As a result, unbounded problems, listed in \Cref{tab::distribution} in row ``unb'', have a significantly higher percentage of expensive \exactsoplex calls. 
On these unbounded problems, the success rate of the ``floating-point'' techniques arithmetic drops significantly to 65.18\% on \fpeasy and 94.89\% on \numdiff.
This leads to an increase in more expensive techniques; most notably, the number of \exactsoplex calls significantly increases from 2.85\% to 29.89\% on \fpeasy and from 0.06\% to 3.86\% on \numdiff.


%% file: sections/conclusion.tex
\section{Conclusion}
\label{sec::conclusion}

Our experiments indicate that the overwhelming majority of the decisions made by the \fpsolver \scip on both the test set \fpeasy and even on the numerically more challenging test set \numdiff are correct even in rational arithmetic.
For feasible instances, only a small, typically single-digit, number of nodes per instance can not be verified and would require further treatment to generate a proof of correctness of the floating-point calculation.
Notably, most of the node decisions can be verified in fast floating-point arithmetic and do not require techniques that rely on more expensive rational arithmetic.

However, we also observe that sub-problems that are infeasible in exact arithmetic pose a special situation.  Here, sometimes the floating-point solver can find solutions with slight violations and use them to prune nodes which are incorrect in terms of rational arithmetic.
This is particularly troubling for instances that are globally infeasible in exact arithmetic.
Our results show that this issue can only partially be addressed by tightening the feasibility tolerance of the solver.


Overall, these results are encouraging.
On the one hand, they quantify the widely accepted belief that floating-point MIP solvers are generally numerically robust for well-behaved input data.
On the other hand, they also suggest a path forward to using floating-point solvers as a grey box in order to solve MIPs exactly over the rational numbers.
%
First, the experimental setup put forward in this paper can be used to generate partial certificates of optimality, which can even be verified for example in the \solver{VIPR} format~\cite{CheungGleixnerSteffy2017}.
Second, these partial certificates can then be completed by continuing to explore sub-problems that were incorrectly discarded, either recursively with increased precision or by calling an exact MIP solver.
For this to become competitive with the state of the art in exact MIP solving, our \emph{a posteriori} verification of LP-based \bandb needs to be extended to include more advanced techniques like presolving or cutting plane separation.

%% file: sections/acknowledgements.tex
\textbf{Acknowledgements.}
The work for this article has been partly conducted within the Research Campus MODAL funded by the German Federal Ministry of Education and Research (BMBF grant number 05M14ZAM).